\documentclass[11pt,fleqn,a4paper]{article}
\usepackage{amssymb,latexsym,amsmath,amsfonts}
\advance\textheight by \topskip

\def\deg{\textrm{deg}}
\def\d{{\rm d}}
\def\e{{\rm e}}

\def\q#1{q_1^{#1}}
\def\p#1{p_1^{#1}}
\def\Q#1{q_2^{#1}}
\def\P#1{p_2^{#1}}
\def\r{r_1}
\def\R{r_2}

\def \N{\mathbb{N}}

\def \C{\mathbb{C}}
\def \O{\mathcal{O}}

\def \reff#1{(\ref{#1})}

\def \f#1#2#3#4#5#6{
      \: {_{#1}} #2_{#3}\left( \left. #4 \atop #5 \right| #6 \right) }

\newtheorem{theorem}{Theorem}[section]
\newtheorem{lemma}[theorem]{Lemma}
\newtheorem{corollary}[theorem]{Corollary}

\newtheorem{Definition}[theorem]{Definition}
\newenvironment{definition}{\begin{Definition}\rm}{\end{Definition}}
\newtheorem{Remark}[theorem]{Remark}
\newenvironment{remark}{\begin{Remark}\rm}{\end{Remark}}
\newtheorem{Example}[theorem]{Example}

\newenvironment{proof}%
{\rm \trivlist \item[\hskip \labelsep{\bf Proof. }]}%
{\hspace*{\fill}$\Box$\endtrivlist}

\newenvironment{varproof}
    {\rm \trivlist \item[\hskip \labelsep{\bf Proof}]}
    {\hspace*{\fill}$\Box$\endtrivlist}

\numberwithin{equation}{section}

\begin{document}

    \begin{center} \Large\bf
        Irrationality proof of certain Lambert series\\ using little
        $q$-Jacobi polynomials
    \end{center}

    \

    \begin{center}
        \large
        J. Coussement,
        C. Smet \\
        \normalsize \em
        Department of Mathematics, Katholieke Universiteit
        Leuven,\\
        Celestijnenlaan 200 B, 3001 Leuven, Belgium \\
        \rm jonathan.coussement@wis.kuleuven.be, christophe.smet@wis.kuleuven.be\\[3ex]
    \end{center}

\ \\[1ex]

\begin{abstract}
    \noindent
We apply the Pad\'{e} technique to find rational approximations to
\[h^{\pm}(q_1,q_2)=\sum_{k=1}^\infty\frac{\q{k}}{1\pm \Q{k}},\qquad
0<q_1,q_2<1,\quad q_1\in\mathbb{Q},\quad q_2=1/p_2,\quad
p_2\in\mathbb{N}\setminus\{1\}.\]
A separate section is dedicated to the special case $q_i=q^{r_i},
r_i\in\mathbb{N},q=1/p,p\in\mathbb{N}\setminus\{1\}$. In this
construction we make use of little $q$-Jacobi polynomials. Our
rational approximations are good enough to prove the irrationality
of $h^{\pm}(q_1,q_2)$ and give an upper bound for the
irrationality measure.
\end{abstract}

\section{Introduction}

In this paper we investigate quantities of the form
\begin{equation}
\label{q-versie} %
h^{\pm}(q_1,q_2):=\sum_{k=1}^\infty\frac{\q{k}}{1\pm \Q{k}},\qquad
0<q_1,q_2<1,\quad q_1\in\mathbb{Q},\quad q_2=1/p_2,\quad
p_2\in\mathbb{N}\setminus\{1\}.
\end{equation}
Since we will assume $q_1 ,q_2$ to be fixed, we will write
$h^{\pm}=h^{\pm}(q_1 ,q_2 )$. In the special case
\begin{equation}
\label{specgeval} %
q_i=q^{r_i},\qquad q=1/p,\quad p\in\mathbb{N}\setminus\{1\},\quad
r_i \in\mathbb{N},\end{equation} by writing $(1+
\Q{k})^{-1}=\sum_{j=0}^\infty (- \Q{k})^j$ and changing the order
of summation, we clearly have
\[\lim_{q\uparrow 1} (1-q)\,h^+=\sum_{j=0}^\infty
\frac{(-1)^j}{\r+j\R}=\frac{1}{\R }\, \Psi(-1,1,\frac{\r }{\R })\]
where $\Psi$ is the Lerch transcendent, which is a generalization
of the Hurwitz zeta function and the polylogarithm function. Some
particular cases are $h^+(q,q)=-\ln_q 2$ and $h^+(q,q^2)=\beta_q
(1)$ which are $q$-extensions of $-\ln 2$ and $\beta(1)=\pi/4$,
respectively. In the same manner $h^-$ can be seen as a
$q$-analogue of the (harmonic) series $\sum_{k=1}^\infty (\r
+k\R)^{-1}$.

In 1948 Erd\H{o}s proved that $h^-(q,q)=\zeta_q(1)$ is irrational
when $q=1/2$, see \cite{Erdos}. Later, Peter Borwein
\cite{bor3,bor1} showed that $\zeta_q(1)$ and $\ln_q 2$ are
irrational whenever $q=1/p$ with $p$ an integer greater than $1$.
Other irrationality proofs were found in, e.g.,
\cite{Amde,Bund,Matala2,Walter1,kelly,Zudilin1,Zudilin2}. To the
best of our knowledge, the sharpest upper bounds for the
irrationality measure of $\zeta_q(1)$ and $\ln_q 2$ which are
known in the literature until now, are $2.42343562$ and
$3.29727451$ respectively \cite{Zudilin1,Zudilin2}.

In \cite{Matala} Matala-aho and Pr\'evost also considered
quantities of the form \reff{q-versie}. However, not all the
numbers we prove to be irrational are covered by their result. To
prove this irrationality we use a well-known lemma, which
expresses the fact that a rational number can be approximated to
order 1 by rational numbers and to no higher order \cite[Theorem
186]{Hardy}.
\begin{lemma}
\label{irrlemma} %
Let x be a real number.  Suppose there exist integers $a_n,b_n
\left(n\in \mathbb{N}\right)$ such that
\begin{enumerate}
    \item[(i)] $b_nx-a_n\neq 0$ for all $n\in \mathbb{N}$;
    \item[(ii)] $\lim\limits_{n\rightarrow\infty}\left(b_nx-a_n\right)=0$,
\end{enumerate}
then x is irrational.
\end{lemma}
\begin{proof}
Suppose $x$ is rational, so write $x=a/b$ with $a,b$ coprime. Then
$b_na-a_nb$ is a nonzero integer sequence that tends to zero,
which is a contradiction.
\end{proof}
In Section~\ref{sectionRA} we construct rational approximations to
$h^{\pm}$. In particular, we extend the Pad\'e approximation
technique applied in \cite{Walter1} to prove the irrationality of
$\zeta_q(1)$ and $\ln_q 2$ and use little $q$-Jacobi polynomials
(which are a generalization of the $q$-Legendre polynomials).
Section~\ref{sectionirr} then mainly consists of calculating the
asymptotic behaviour of the 'error term'.
Section~\ref{sectionspec} points out what improvements can be made
in the special case \reff{specgeval}. If we define
\begin{equation}
\label{eta}\eta^-:=1+\frac{3}{\pi^2},\qquad
\eta^+:=1+\frac{4}{\pi^2},
\end{equation}
\begin{equation}
\label{gamma-}
\gamma^{-}(\R ):=\frac{3}{\pi^2}\left[1+2\prod_{\substack{\varpi|\R\\
\varpi\: {\rm
prime}}}\frac{\varpi^2}{\varpi^2-1}\sum_{\substack{l=1\\ (l,\R
)=1}}^{\R }\frac{1}{l^2}-\frac{1}{\R}\prod_{\substack{\varpi|\R\\
\varpi\: {\rm
prime}}}\frac{\varpi}{\varpi+1}\right]<1+\frac{3}{\pi^2}
\end{equation}
and
\begin{multline}
\label{gamma+} %
\gamma^+(\R ):=\frac{1}{\pi^2}\left[4
+6\prod_{\substack{\varpi|\R\\
\varpi\: {\rm
prime}}}\frac{\varpi^2}{\varpi^2-1}\sum_{\substack{l=1\\ (l,\R
)=1}}^{\R }\frac{1}{l^2} \right] \\
-\frac{1-(-1)^{\R}}{2\pi^2}\left[ \frac{2}{\R}\prod_{\substack{\varpi|\R\\
\varpi\: {\rm prime}}}\frac{\varpi}{\varpi+1}
+\prod_{\substack{\varpi|\R\\
\varpi\: {\rm prime}}}\frac{\varpi^2}{\varpi^2-1}
\sum_{\substack{l=\lceil \frac{\R{}}{2} \rceil \\(l,\R )=1}}^{\R }
\frac{1}{l^2}\right] < 1+\frac{4}{\pi^2},
\end{multline}
then our main results are the following.
\begin{theorem}
\label{irrationality alg}%
Let $q_2=1/p_2$ with $p_2\in\mathbb{N}\setminus\{1\}$ and
$q_1\in\mathbb{Q}$ with $0<q_1<1$. Then the number $h^{\pm}$,
defined as in \reff{q-versie}, is irrational. Moreover, there
exist integer sequences $a_{n}^\pm$, $b_{n}^\pm$ such that
\begin{equation}
\label{rest asymp in p alg} %
\lim_{n\rightarrow\infty}\left|b_n^\pm h^{\pm}-a_n^\pm
\right|^{1/n^2}\leq \P{\eta^\pm-\frac{3}{2}}<1
\end{equation}
and
\begin{equation}
\label{asympbeta alg}%
 \lim_{n\rightarrow\infty}|b_n^\pm|^{1/n^2}\leq
\P{\eta^\pm+\frac{3}{2}}.
\end{equation}
\end{theorem}
\begin{theorem}
\label{irrationality}%
Let $q=1/p$ with $p\in\mathbb{N}\setminus\{1\}$, $q_i^{}=q^{r_i}$
and $p_i^{}=p^{r_i}$ with $r_i\in\mathbb{N}$, $i=1,2$ and
$(r_1,r_2)=1$. Then the number $h^{\pm}$, defined as in
\reff{q-versie}, is irrational. Moreover, there exist integer
sequences $a_{n}^\pm$, $b_{n}^\pm$ such that
\begin{equation}
\label{rest asymp in p} %
\lim_{n\rightarrow\infty}\left|b_n^\pm h^{\pm}-a_n^\pm
\right|^{1/n^2}\leq \P{\gamma^\pm(\R )-\frac{3}{2}}<1
\end{equation}
and
\begin{equation}
\label{asympbeta}%
 \lim_{n\rightarrow\infty}|b_n^\pm|^{1/n^2}\leq
\P{\gamma^\pm(\R )+\frac{3}{2}}.
\end{equation}
\end{theorem}
\begin{remark}
In fact, Theorem~\ref{irrationality} can also be applied if ${\rm
gcd}(\r,\R)=\rho\not=1$. In that case just note that
$h^{\pm}(q_1,q_2)=h^{\pm}(q'^{r_1'},q'^{r_2'})$ with
$\r'=\r/\rho$, $\R'=\R/\rho$ and $q'=q^\rho$.
\end{remark}

As a side result of the irrationality, we also obtain an upper
bound for the irrationality measure (Liouville-Roth number, order
of approximation) for $h^{\pm}$. Recall that this measure is
defined as
\[\mu(x):=\inf\left\{t\ :\ \left|x-\frac{a}{b}\right|>\frac{1}{b^{t+\varepsilon}},\
\forall\varepsilon>0,\ \forall a,b \in\mathbb{Z},\ b
\:\rm{sufficiently}\: \rm{large}\right\},\]
see, e.g., \cite{bor2}. It is known that all rational numbers have
irrationality measure 1, whereas irrational numbers have
irrationality measure at least 2. Furthermore, if $b_nx-a_n\neq 0$
for all $n\in \mathbb{N}$, $\left|b_nx-a_n\right|=\O(b_n^{-s})$
with $0<s<1$  and $|b_n|<|b_{n+1}|<|b_n|^{1+o(1)}$, then the
measure of irrationality satisfies $2\le\mu(x)\le 1+1/s$, see
\cite[exercise 3, p.~376]{bor2}. Note that by \reff{rest asymp in
p alg} and \reff{asympbeta alg}, respectively \reff{rest asymp in
p} and \reff{asympbeta}, we get the asymptotic behaviour
\begin{equation}
\left|b_n^\pm h^{\pm}-a_n^\pm
\right|=\O\Bigl((b_n^\pm)^{-\frac{3-2\,\eta^\pm}{3+2\,\eta^\pm
}+\varepsilon}\Bigr),\qquad \mbox{for all } \varepsilon
>0,\qquad n\to \infty,
\end{equation}
and for the special case \reff{specgeval}
\begin{equation}
\left|b_n^\pm h^{\pm}-a_n^\pm
\right|=\O\Bigl((b_n^\pm)^{-\frac{3-2\,\gamma^\pm (\R
)}{3+2\,\gamma^\pm (\R )}+\varepsilon}\Bigr),\qquad \mbox{for all
} \varepsilon
>0,\qquad n\to \infty,
\end{equation}
which then implies the following upper bound for $\mu(h^{\pm})$.
\begin{corollary}
\label{irrationality measure} %
Under the assumptions of Theorem~\ref{irrationality alg} we have
$2\le \mu(h^{\pm})\le \nu^\pm$; under the assumptions of
Theorem~\ref{irrationality} we have $2\le \mu(h^{\pm})\le m^\pm
(\R)$ where
\begin{equation}
\label{upb} %
m^\pm (\R) = \left(\frac{3-2\,\gamma^\pm (\R )}{6}\right)^{-1},
\end{equation}
with $m^\pm (\R)\le \nu^\pm$, where $\nu^+=\frac{6\pi^2}{\pi^2-8}$
and $\nu^-=\frac{6\pi^2}{\pi^2-6}$.

\end{corollary}
\begin{table}[t]
\begin{center}
\begin{tabular}{|c|rcl|rcl|}
\hline %
$\R$ &  \multicolumn{3}{|c|}{$m^-(\R)$}&\multicolumn{3}{|c|}{$m^+(\R)$}\\
\hline & & & & & & \\[-2ex]
1 & $\frac{2\pi^2}{\pi^2-4}$&$\approx$&$ 3.362953864 $&$\frac{6\pi^2}{3\pi^2-14}$&$\approx$&$ 3.793858357 $ \\[0.5ex]
2 & $\frac{6\pi^2}{3\pi^2-20}$&$\approx$&$ 6.162845000 $&$\frac{2\pi^2}{\pi^2-8}$&$\approx$&$ 10.55796017 $ \\[0.5ex]
3 & $\frac{16\pi^2}{8\pi^2-57}$&$\approx$&$ 7.192005083 $&$\frac{96\pi^2}{48\pi^2-373}$&$\approx$&$ 9.405127174$ \\[0.5ex]
4 & $\frac{54\pi^2}{27\pi^2-205}$&$\approx$& $8.668909282 $&$\frac{54\pi^2}{27\pi^2-232}$&$\approx$&$ 15.45734242 $ \\[0.5ex]
5 & $\frac{1728\pi^2}{864\pi^2-6565}$&$\approx$&$ 8.690997496 $&$\frac{10368\pi^2}{5184\pi^2-42797}$&$\approx$&$ 12.22991528 $ \\[0.5ex]
6 & $\frac{300\pi^2}{150\pi^2-1211}$&$\approx$& $10.98899223 $&$\frac{150\pi^2}{75\pi^2-668}$&$\approx$&$ 20.49894619 $ \\[0.5ex]
7 & $\frac{86400\pi^2}{43200\pi^2-338681}$&$\approx$&$ 9.724867074 $&$\frac{103680\pi^2}{51840\pi^2-440701}$&$\approx$&$ 14.42473632 $ \\[0.5ex]
8 & $\frac{132300\pi^2}{66150\pi^2-534587}$&$\approx$& $11.03878708 $&$\frac{66150\pi^2}{33075\pi^2-294856}$&$\approx$&$ 20.67290169 $ \\[0.5ex]
9 & $\frac{940800\pi^2}{470400\pi^2-3801647}$&$\approx$&$ 11.04061736 $&$\frac{1128960\pi^2}{564480\pi^2-4937467}$&$\approx$&$ 17.58230823 $ \\[0.5ex]
10 & $\frac{71442\pi^2}{35721\pi^2-294473}$&$\approx$& $12.14040518 $&$\frac{71442\pi^2}{35721\pi^2-322256}$&$\approx$&$ 23.27373406 $ \\[0.5ex]
\hline%
\end{tabular}
\caption{\label{tabel1} Some values of the upper bound $m^\pm(\R)$
for the irrationality measure of $h^\pm$.}
\end{center}
\end{table}
\begin{remark}
In the case \reff{specgeval} with $\R =1$, we can sharpen the
upper bound $m^\pm(\R)$. We will discuss this in
Section~\ref{finalremark}. In particular, we will show that
$\mu(\zeta_q(1))\le \frac{2\pi^2}{\pi^2-2}\approx 2.508284762$,
which was also found in \cite{Walter1}, and $\mu(\ln_q 2)\le
\frac{6\pi^2}{3\pi^2-8}\approx 2.740438628$, which is a better
upper bound than the one in \cite{Zudilin1}.
\end{remark}

\section{Rational approximation}
\label{sectionRA}

We first focus on the general case \reff{q-versie}, the special
case \reff{specgeval} will be treated in
Section~\ref{sectionspec}.  We use the notation $q_1=s_1/t_1$ with
$\gcd(s_1,t_1)=1$, and $p_1=1/q_1$.

\subsection{Pad\'e approximation}

To prove the irrationality of $h^{\pm}$ we will apply
Lemma~\ref{irrlemma}. So we need a sequence of 'good' rational
approximations. To find these we will perform the (well-known)
idea of Pad\'e approximation to the Markov function
\begin{equation}
\label{Markov function} f(z):=\sum_{k=0}^\infty
\frac{\q{k}}{z-\Q{k}}
=\int_0^1\frac{q_1^{\log_{q_2}x}}{z-x}\, \frac{\d_{\Q{}}x}{x}, %
\end{equation}
where $\log_qx=\frac{\log x}{\log q}$ and the $q$-integration is
defined as
\begin{equation}
\int_0^1 g(x)\,\d_q x:=\sum_{k=0}^\infty q^k g(q^k).
\end{equation}
So, we look for polynomials $P_n$ and $Q_n$ of degree $n$ such
that
\begin{equation}
\label{PA} %
Q_n(z)f(z)-P_n(z)=O\left(z^{-n-1}\right), \qquad
z\rightarrow\infty.
\end{equation}
As is well known in the Pad\'e approximation theory (see, e.g.,
\cite{Nikishin}) the polynomials $Q_n$ then satisfy the
orthogonality relations
\begin{equation}
\label{OP} %
\int_0^1 Q_n(x)\, x^m \, q_1^{\log_{q_2}x}\,x^{-1}\,
\d_{\Q{}}x=0,\qquad m=0,\dots,n-1.
\end{equation}
Little $q$-Jacobi polynomials satisfy
\begin{equation}
\label{defortholittleqjac} %
\sum_{k=0}^\infty
p_n(q^k;a,b|q)p_m(q^k;a,b|q)(aq)^k\frac{(bq;q)_k}{(q;q)_k}=0,\qquad
m\neq n.
\end{equation}
Hence $Q_n$ are little $q$-Jacobi polynomials with a particular
set of parameters, namely $Q_n(z)=p_n(z;\q{}\P{},1|\Q{})$, see,
e.g., \cite[Section~3.12]{Koekoek}.
\begin{lemma}
The polynomials $Q_n$ have the explicit expressions
\begin{eqnarray}
\label{polynomialQ1} %
Q_n(z) & = & %
\sum_{k=0}^n \frac{(\P{n};\Q{})_k (\q{}
\Q{n};\Q{})_k}{(\q{};\Q{})_{k}(\Q{};\Q{})_k}\,\Q{k}\,
z^k,\\
\label{polynomialQ2} %
& = & \frac{(\P{n};\Q{})_n}{(\q{};\Q{})_n}\sum_{k=0}^n
\frac{(\P{n};\Q{})_k(\q{}\Q{n};\Q{})_k}{\left[(\Q{};\Q{})_{k}\right]^2}\,\Q{k}\,
(\Q{}z;\Q{})_k.
\end{eqnarray}
\end{lemma}
\begin{proof}
From, e.g., \cite[Section 3.12]{Koekoek}, we know that the
polynomials satisfying the orthogonality conditions \reff{OP} have
the hypergeometric expression
\[
Q_n(z)=\f{2}{\phi}{1}{\P{n},\q{}\Q{n}}{\q{}}{\Q{};\Q{}z},
\]
which is \reff{polynomialQ1}. Next we apply the transformation
formula \cite[(0.6.24)]{Koekoek} and find
\[
Q_n(z)=\frac{(\P{n};\Q{})_n}{(\q{};\Q{})_n}
\f{3}{\phi}{2}{\P{n},\q{}\Q{n},\Q{}z}{\Q{},0}{\Q{};\Q{}},
\]
giving the expression \reff{polynomialQ2}.
\end{proof}
It is easily checked that the $P_n$ are connected with the
polynomials $Q_n$ by the formula
\begin{equation}
\label{pade geeft Pn} %
P_n(z)=\int_0^1\frac{Q_n(z)-Q_n(x)}{z-x}\,
q_1^{\log_{q_2}x}\,x^{-1}\, \d_{\Q{}}x.
\end{equation}
Indeed, then
\begin{equation}
\label{errorint} %
Q_n(z)f(z)-P_n(z)=\int_0^1\frac{Q_n(x)}{z-x}\,
q_1^{\log_{q_2}x}\,x^{-1}\, \d_{\Q{}}x,
\end{equation}
and the conditions \reff{PA} are fulfilled.
\begin{lemma}
The polynomials $P_n$ have the explicit formulae
\begin{eqnarray}
\label{polynomialP1}%
P_n(z) & = & \sum_{k=0}^n\frac{(\P{n};\Q{})_k
(\q{}\Q{n};\Q{})_k}{(\q{};\Q{})_{k}(\Q{};\Q{})_k}\,\Q{k}\,
\sum_{j=0}^{k-1}\frac{z^j}{1-\q{}\Q{k-j-1}},\\[1ex]
\nonumber %
& = & -\frac{(\P{n};\Q{})_n}{(\q{};\Q{})_n}\sum_{k=0}^n
\frac{(\P{n};\Q{})_k (\q{}\Q{n};\Q{})_k}{\left[(\Q{};\Q{})_{k}\right]^2}\,\Q{k}\\
\label{polynomialP2} %
& & \hspace{4cm} \times
\sum_{j=1}^k\Q{j}\,(\Q{j+1}z;\Q{})_{k-j}\,\frac{(\Q{};\Q{})_{j-1}}{(\q{};\Q{})_j}.
\end{eqnarray}
\end{lemma}
\begin{proof}
First of all we mention that (for $k\in \N\cup \{0\}$)
\[\int_0^1\frac{z^k-x^k}{z-x}\, q_1^{\log_{q_2}x}\,x^{-1}\,
\d_{\Q{}}x=
\sum_{j=0}^{k-1}z^j\sum_{\ell=0}^{\infty}\Q{(k-1-j)\ell}\q{\ell}=\sum_{j=0}^{k-1}\frac{z^j}{1-\q{}\Q{k-j-1}}.\]
So, applying \reff{polynomialQ1} to \reff{pade geeft Pn} we easily
obtain \reff{polynomialP1}. Next, observe that
\begin{equation}
\label{z-x} %
\frac{(\Q{} z;\Q{})_k-(\Q{}
x;\Q{})_k}{z-x}=-\sum_{j=1}^k\Q{j}\,(\Q{j+1}z;\Q{})_{k-j}\,(\Q{}
x;\Q{})_{j-1}
\end{equation}
which one can prove by induction. Moreover, using the $q$-binomial
series \cite[Section~10.2]{Andrews}, \cite[Section~1.3]{Gasper}
\[\sum_{n=0}^\infty\frac{(a;q)_n}{(q;q)_n}x^n=\frac{(ax;q)_\infty}{(x;q)_\infty},\qquad |q|<1,\ |x|<1,\]
we get
\begin{equation}
\label{integral J} %
\int_0^1(\Q{} x;\Q{})_{j-1}\, q_1^{\log_{q_2}x}\,x^{-1}\,
\d_{\Q{}}x= (\Q{};\Q{})_{j-1}\sum_{\ell=0}^\infty
\frac{(\Q{j};\Q{})_{\ell}}{(\Q{};\Q{})_{\ell}}\,\q{\ell} =
\frac{(\Q{};\Q{})_{j-1}}{(\q{};\Q{})_j}.
\end{equation}
Combining \reff{pade geeft Pn}, \reff{polynomialQ2}, \reff{z-x}
and \reff{integral J} we then finally establish
\reff{polynomialP2}.
\end{proof}

\subsection{Rational approximants to $h^{\pm}$}
\label{an bn integer}

Notice that by the definition \reff{Markov function} of $f$ we
have
\[h^{\pm}=\mp \frac{\q{}}{\Q{}}\, f(\mp \P{}).\]
Following the idea of Pad\'e approximation we could try to
approximate $h^{\pm}$ by the sequence of rational numbers $\mp
\q{} P_n(\mp \P{})/[\Q{}Q_n(\mp \P{})]$. However, we prefer the
evaluation of $f$ at $\mp \P{n}$, which gives
\begin{equation}
\label{f in -p2n}%
h^{\pm}=\sum_{k=1}^{n-1}\frac{\q{k}}{1\pm \Q{k}} \mp
\left(\frac{\q{}}{\Q{}}\right)^n f(\mp \P{n}).
\end{equation}
In this way we can benefit from the fact that the finite sum on
the right hand side of \reff{f in -p2n} already gives a good
approximation for $h^{\pm}$. Moreover, the approximation \reff{PA}
is useful for $z$ tending to infinity. Hence the evaluation at the
point $\mp \P{n}$ makes more sense, especially since $n$ itself
will tend to infinity too.
So, a 'natural' choice for rational approximations
$a_n^{\pm}/b_n^{\pm}$ to $h^{\pm}$ then has the following
expressions.
\begin{definition}
\label{defanbn} %
Define
\begin{align}
\label{an} & a_n^\pm :=  e_n^\pm\left[\p{n}\, Q_n(\mp
\P{n})\sum_{k=1}^{n-1}\frac{\q{k}}{1\pm \Q{k}}\mp \P{n}\,P_n(\mp
\P{n})\right],\\[1ex]
\label{bn} & b_n^\pm :=  e_n^\pm\ \p{n}\,Q_n(\mp \P{n}),
\end{align}
where $e_n^\pm$ are factors such that these are integer sequences.
\end{definition}
The following lemma gives a possible choice for the factors
$e_n^\pm$.
\begin{lemma}
\label{lemma en alg} %
By taking
\begin{equation}
\label{en alg} %
e_n^\pm= %
\left(\prod_{k=0}^{n-1}(t_1p_2^k-s_1)\right)^2\,s_1^n\,{\rm
lcm}\left\{\left. p_2^k\pm1 \, \right| \, 0\le k \le n-1 \right\}
\end{equation}
the $a_n^\pm$ and $b_n^\pm$, defined as in \reff{an} and
\reff{bn}, are integer sequences.
\end{lemma}
\begin{proof}
It is not very convenient to prove that an expression is an
integer when it depends on the rational $q_2$. So, we first write
$Q_n(\mp \P{n})$ depending on the integer $\P{}$. By
\reff{polynomialQ1} we get
\begin{equation}
\label{Qn in p} %
Q_n(\mp \P{n})  =  \sum_{k=0}^n \frac{(\P{n};\Q{})_k
(\p{}\P{n};\P{})_k}{(\p{};\P{})_{k}(\P{};\P{})_k}
\,\P{\frac{k^2-k}{2}}(\pm 1)^k.
\end{equation}
Notice that
\[\frac{(\P{n};\Q{})_k}{(\P{},\P{})_k} = \frac{(\P{};\P{})_n}{(\P{};\P{})_{n-k}(\P{};\P{})_k}
=\left[{n\atop k}\right]_{\P{}},\]
which is an integer. Moreover, the possible denominators $t_1$
appear as often in the numerator as in the denominator of $Q_n(\mp
\P{n})$.  The factor $s_1^n$ in $e_n^\pm$ is needed because of the
factor $p_1^n$ in $a_n^\pm$ and $b_n^\pm$.  So the only
denominators in $Q_n(\mp \P{n})$ originate from $(\p{};\P{})_{k}$,
and hence they are cancelled out by the first product in
$e_n^\pm$. This already implies that $b_n^\pm$ is an integer.
Obviously, then also
\[e_n^\pm\,\p{n}\,
Q_n(\mp \P{n})\sum_{k=1}^{n-1}\frac{\q{k}}{1\pm \Q{k}}=b_n^\pm
\sum_{k=1}^{n-1}\frac{\q{k}}{1\pm \Q{k}}\]
is an integer by the definition of $e_n^\pm$. So, what remains to
prove is that $e_n^\pm\, \P{n}\,P_n(\mp \P{n})$ is an integer. By
\reff{polynomialP1} we have
\begin{equation}
\label{Pn een}%
P_n(\mp \P{n})=\sum_{k=0}^n\frac{(\P{n};\Q{})_k
(\p{}\P{n};\P{})_k}{(\p{};\P{})_{k}(\P{};\P{})_k}\,(-1)^k\,
\sum_{j=0}^{k-1}\frac{(\mp 1)^{j}\,
\p{}\P{\frac{k^2+k}{2}+n(j-k)-j-1}}{\p{}\P{k-j-1}-1}.
\end{equation}
Since $(\P{};\P{})_k$ is a divisor of $(\P{n};\Q{})_k$, it is
clear that by \reff{en} the only possible denominators in $e_n^\pm
\, \P{n}\,P_n(\mp \P{n})$ are powers of $p_2$. The formula
\reff{polynomialP2} leads to
\begin{multline*}
P_n(\mp \P{n})=\frac{(\P{n};\Q{})_n}{(\p{};\P{})_n}\sum_{k=0}^n
\frac{(\P{n};\Q{})_k
(\p{}\P{n};\P{})_k}{\left[(\P{};\P{})_{k}\right]^2}\,\p{n-k}\P{\binom{n-k}{2}}\,(-1)^{n-k}\\
\times
\sum_{j=1}^k\left(\frac{\p{}}{\P{}}\right)^{j}\,(\mp\P{n-j-1};\Q{})_{k-j}\,\frac{(\P{};\P{})_{j-1}}{(\p{};\P{})_j}.
\end{multline*}
\end{proof}

\begin{remark}
Looking at \reff{Pn een} one could expect a power of $p_2$ in the
denominator of $\P{n}\,P_n(\mp \P{n})$, and this of the order of
$\P{n^2/2}$. This would totally ruin the asymptotics in the next
section. However, Maple calculations showed the absence of a power
of $p_2$ in the denominator. This is why we had to use an
equivalent formula for $P_n$, which is given by
\reff{polynomialP2}.
\end{remark}

\section{Irrationality of $h^{\pm}$}
\label{sectionirr}

In this section we look at the error term $b_n^\pm h^\pm -
a_n^\pm$, where $a_n^\pm$ and $b_n^\pm$ are defined as in
Definition~\ref{defanbn} and \reff{en alg}. Using \reff{f in -p2n}
and \reff{errorint} one easily sees that it has the integral
representation
\begin{align}
\nonumber %
b_n^\pm h^\pm - a_n^\pm & = \mp e_n^\pm\, \P{n} \Bigl[Q_n(\mp
\P{n})f(\mp \P{n})-P_n(\mp \P{n})\Bigr]\\[1ex]
& = e_n^\pm \, \P{n} \int_0^1\frac{Q_n(x)}{\P{n}\pm x}\,
q_1^{\log_{q_2}x}\,x^{-1}\, \d_{\Q{}}x.
\end{align}
We will show that this expression is different from zero for all
$n\in \N$ and obtain its asymptotic behaviour. Here we study
\begin{equation}
\label{Rn} %
R_n^\pm:=\int_0^1\frac{Q_n(x)}{\P{n}\pm x}\,
q_1^{\log_{q_2}x}\,x^{-1}\, \d_{\Q{}}x
\end{equation}
and $e_n^\pm$ separately.

\subsection{Asymptotic behaviour of $R_n^\pm$} %
\label{restterm asymp}
We will need the following very general lemma for sequences of
polynomials with uniformly bounded zeros. This can be found in,
e.g., \cite[Lemma 3]{Walter1}, but we include a short proof for
completeness.
\begin{lemma}
\label{hulplemma} %
Let $\{\pi_n\}_{n\in \N}$ be a sequence of monic polynomials for
which $\deg(\pi_n)=n$ and the zeros $x_{j,n}$ satisfy
$|x_{j,n}|\le M$, with $M$ independent of $n$. Then
\[\lim_{n\to \infty} \left|\pi_n\left(cx^n\right)\right|^{1/n^2}=|x|,
\qquad |x|>1,\ c\in\C.\]
\end{lemma}
\begin{proof}
Since $|x|>1$, for large $n$ we easily get
\[0\le |cx^n|-M\le |cx^n-x_{j,n}|\le |cx^n|+M, \qquad
j=1,\ldots,n.\]
This implies
\[\left(|cx^n|-M\right)^n\le
\left|\pi_n\left(cx^n\right)\right|\le \left(|cx^n|+M\right)^n\]
and
\[|x|\left(|c|-\frac{M}{|x|^n}\right)^{1/n}\le
\left|\pi_n\left(cx^n\right)\right|^{1/n^2}\le
|x|\left(|c|+\frac{M}{|x|^n}\right)^{1/n}.\]
The lemma then follows by taking limits.
\end{proof}
For $R_n^\pm$, defined as in \reff{Rn}, we have the following
asymptotic result. Here we use a similar reasoning as in
\cite{Walter1} for the irrationality of $\zeta_q(1)$.
\begin{lemma}
\label{restterm asymp lemma} %
Let $Q_n$ be the polynomials \reff{polynomialQ1} satisfying the
orthogonality relations \reff{OP}. Then $R_n^\pm$ is different
from zero for all $n$ and
\begin{equation}
\label{Rnasymp} \lim\limits_{n\to \infty}
\left|R_n^\pm\right|^{1/n^2}=\P{-3/2}.
\end{equation}
\end{lemma}
\begin{proof}
First of all observe that
\[
Q_n(\mp \P{n})\,R_n^\pm = \mp \int_0^1 Q_n(x)\,\frac{Q_n(\mp
\P{n})-Q_n(x)}{\mp \P{n}-x} \, q_1^{\log_{q_2}x}\,x^{-1}\,
\d_{\Q{}} x + \int_0^1 \frac{Q_n^2(x)}{\P{n}\pm x}\,
q_1^{\log_{q_2}x}\,x^{-1}\, \d_{\Q{}} x.
\]
The first integral on the right hand side vanishes because of the
orthogonality relations \reff{OP} for the polynomial $Q_n$.
Furthermore, note that $0\le x\le 1$ so that
\begin{equation}
\label{afschatting} %
0<\frac{1}{\P{n}+1}\int_0^1 Q_n^2(x)\, q_1^{\log_{q_2}x}\,x^{-1}\,
\d_{\Q{}} x \le Q_n(\mp \P{n})\,R_n^\pm \le
\frac{1}{\P{n}-1}\int_0^1 Q_n^2(x)\, q_1^{\log_{q_2}x}\,x^{-1}\,
\d_{\Q{}} x.
\end{equation}
This already proves that $R_n^\pm\not=0$. Next, from
\cite[(3.12.2)]{Koekoek} we get
\[\int_0^1 Q_n^2(x)\, q_1^{\log_{q_2}x}\,x^{-1}\,
\d_{\Q{}} x=
\frac{\q{n}}{1-\q{}\Q{2n}}\left(\frac{(\Q{};\Q{})_n}{(\q{};\Q{})_n}\right)^2.\]
Applying this on \reff{afschatting}, we easily establish
\begin{equation}
\label{R naar Q} %
\lim\limits_{n\to \infty}\left|Q_n(\mp
\P{n})\,R_n^\pm\right|^{1/n^2}=1.
\end{equation}
Now write $Q_n(x)=\kappa_n\, \hat Q_n(x)$ where $\hat Q_n$ is
monic. From \reff{polynomialQ1} we get that the leading
coefficient $\kappa_n$ has the expression
\[\kappa_n=\frac{(\P{n};\Q{})_n(\q{}\Q{n};\Q{})_n}{(\q{};\Q{})_{n}(\Q{};\Q{})_{n}}\,\Q{n}.\]
Since $\prod_{i=1}^n \P{i-1}\le |(\P{n};\Q{})_n| \le \prod_{i=1}^n
\P{i}$, this gives the asymptotic behaviour
\begin{equation}
\label{kappa} %
\lim_{n\to \infty}|\kappa_n|^{1/n^2}=\P{1/2}.
\end{equation}
Since the $Q_n$ are orthogonal polynomials with respect to a
positive measure on $[0,1]$, their zeros are all in $[0,1]$. From
Lemma~\ref{hulplemma} we then also get
\begin{equation}
\label{hat Q} %
\lim\limits_{n\to \infty}\left|\hat Q_n\left(\mp
\P{n}\right)\right|^{1/n^2}=\P{}.
\end{equation}
Applying \reff{kappa} and \reff{hat Q} to \reff{R naar Q} then
completes the proof.
\end{proof}

\subsection{Asymptotic behaviour of $e_n^\pm$}

We obviously have the asymptotic properties

\begin{align}
\label{enasymp alg 1} %
&
\lim_{n\rightarrow\infty}\left(\prod_{k=0}^{n-1}(t_1p_2^k-s_1)\right)^{1/n^2}=p_2^{1/2}
,\\[1ex]
\label{enasymp alg 2} %
& \lim_{n\rightarrow\infty}\left(s_1^n\right)^{1/n^2}=1
,\\[1ex]
\label{enasymp alg 3} %
& \lim_{n\rightarrow\infty}\left({\rm lcm}\left\{\left. p_2^k-1 \,
\right| \, 0\le k \le n-1 \right\}\right)^{1/n^2}\leq p_2^{3/\pi^2} ,\\[1ex]
\label{enasymp alg 4} %
& \lim_{n\rightarrow\infty}\left({\rm lcm}\left\{\left. p_2^k+1 \,
\right| \, 0\le k \le n-1 \right\}\right)^{1/n^2}\leq
p_2^{4/\pi^2},
\end{align}
where the latter two are well-known properties of the least common
multiple that can easily be deduced from the asymptotic results as
given in \reff{Azonderb}. This leads us to the following
asymptotic behaviour of $e_n^\pm$.

\begin{corollary}
\label{en asymp alg} %
For $e_n^\pm$ defined as in \reff{en alg} we have
\begin{equation}
\label{en asymp formule alg} %
\lim_{n\rightarrow\infty}\left|e_n^\pm\right|^{1/n^2} \le
\P{\eta^\pm},
\end{equation}
where $\eta^\pm$ is defined as in \reff{eta}.  As a result we now
have
\begin{equation}
\label{restterm asymp lemma nieuw}
\lim_{n\rightarrow\infty}\left|b_n^\pm h^{\pm}-a_n^\pm
\right|^{1/n^2}\le \P{\eta^\pm-\frac{3}{2}}.
\end{equation}
\end{corollary}

\subsection{Proof of Theorem~\ref{irrationality alg}}

In the previous sections we defined integer sequences $a_n^\pm$
and $b_n^\pm$ and managed to find the asymptotic behaviour of
$b_n^\pm h^{\pm}-a_n^\pm$.  Putting these results together, we can
now prove Theorem~\ref{irrationality alg}.
\begin{varproof}{\bf of Theorem~\ref{irrationality alg}.}
In Lemma~\ref{lemma en alg} we made sure that $a_n^\pm$ and
$b_n^\pm$, defined as in \reff{an} and \reff{bn}, are integer
sequences. Note that by \reff{kappa}, \reff{hat Q} and \reff{en
asymp formule alg} we then get
\begin{equation}
\lim_{n\rightarrow\infty} |b_n|^{1/n^2}\leq
\P{\eta^\pm+\frac{3}{2}}.
\end{equation}
Lemma~\ref{restterm asymp lemma} assures us that $b_n^\pm
h^{\pm}-a_n^\pm\neq 0$ for all $n\in \mathbb{N}$ and since
$\eta^\pm<\frac32$, \reff{restterm asymp lemma nieuw} guarantees
that $\lim_{n\rightarrow\infty}\left|b_n^\pm
h^{\pm}-a_n^\pm\right|=0$. So, all the conditions of
Lemma~\ref{irrlemma} are fulfilled and $h^{\pm}$ is irrational.
\end{varproof}

\section{Improvements on the results in the special case \reff{specgeval}}
\label{sectionspec}

Throughout this section we consider the special case given by
\reff{specgeval}.  The only difference with the general case is
that we can (in some cases considerably) improve the factor
$e_n^\pm$, which is needed to make the approximation sequences
into integer sequences. The following lemma gives the enhanced
formula for $e_n^\pm$, and can be seen as an analogue of
Lemma~\ref{en alg}.

\begin{lemma}
\label{lemma en} %
By taking
\begin{multline}
\label{en} %
e_n^\pm= %
{\rm lcm}\left\{\left. {\rm denom}
\left((\p{}\P{n};\P{})_k\,(\p{};\P{})_{k}^{-1}\right) \, \right|
\, 0\le k \le n-1 \right\}\\[1ex]
\times {\rm lcm}\left\{\left. \P{j}\pm 1, \, \p{}\P{k}-1 \,
\right| \, 1\le j\le n-1, \, 0\le k \le n-1 \right\}
\end{multline}
the $a_n^\pm$ and $b_n^\pm$, defined as in \reff{an} and
\reff{bn}, are integer sequences.
\end{lemma}
\begin{proof}
The proof is completely analogous to the proof of Lemma~\ref{lemma
en alg}.  There is no factor $s_1^n$ needed since in this case
$p_1$ is an integer.

\end{proof}

\subsection{Asymptotic behaviour of $e_n^\pm$}

In order to obtain some asymptotic results for the quantities
$e_n^\pm$, see \reff{en}, we will use the cyclotomic polynomials
\begin{equation}
\label{cyclotomic def} %
\Phi_n(x)=\prod_{\substack{k=1,\\(k,n)=1}}^n\left(x-e^{\frac{2\pi
ik}{n}}\right).
\end{equation}
Their degree is denoted by Euler's totient function $\phi(n)$,
being the number of positive integers $\le n$ that are coprime
with $n$. It is well-known \cite[Section~4.8]{Still} that
\begin{equation}
\label{nice P} %
x^n-1=\prod_{d|n}\Phi_d(x),\qquad  n=\sum_{d|n}\phi (d),%
\end{equation}
and that every cyclotomic polynomial is monic, has integer
coefficients and is irreducible over $\mathbb{Q} [x]$.
Furthermore, some interesting asymtotic properties are
\begin{align}
\label{Azonderb} %
& \lim_{n\to \infty}\frac{1}{n^2}
\sum_{j=0}^n\phi(aj) =\frac{3a}{\pi^2}\prod_{\substack{\varpi|a\\
\varpi\: {\rm
prime}}}\frac{\varpi}{\varpi+1},\\[1ex]
\label{Ametb} %
& \lim_{n\to \infty}\frac{1}{n^2} \sum_{j=0}^n\phi(aj+b)
=\frac{3a}{\pi^2}\prod_{\substack{\varpi
|a\\\varpi \: {\rm prime}}}\frac{\varpi^2}{\varpi^2-1},\\[1ex]
\label{A2metb} %
& \lim_{n\to \infty}\frac{1}{n^2} \sum_{j=0}^n\phi(2(aj+b))
=\frac{4a}{\pi^2}\prod_{\substack{\varpi |a,\, \varpi\ge 3\\\varpi
\: {\rm prime}}}\frac{\varpi^2}{\varpi^2-1},
\end{align}
where $(a,b)=1$, see~\cite{Bavencoffe,Bezivin,Matala}. They imply
the following results.
\begin{lemma}
\label{lemmaLCM-} %
Let $\r ,\R \in\mathbb{N}$ and $(\r ,\R )=1$. Then
\begin{equation}
\label{AsLCM-} %
\lim_{n\to \infty}\left[{\rm lcm}\left\{\left. \P{j}- 1, \,
\p{}\P{k}-1 \, \right| \, 1\le j\le n-1, \, 0\le k \le n-1
\right\}\right]^{1/n^2}\le \P{\theta^{-} (\R )},
\end{equation}
where
\[
\theta^-(\R )=\gamma^-(\R ) - \frac{3}{\pi^2}\prod_{\substack{\varpi|\R\\
\varpi\: {\rm
prime}}}\frac{\varpi^2}{\varpi^2-1}\sum_{\substack{l=1\\ (l,\R
)=1}}^{\R }\frac{1}{l^2}
\]
with $0<\theta^{-}(\R )<\frac{3}{\pi^2}+\frac{1}{2}$ and
$\gamma^-(\R )$ defined as in \reff{gamma-}.
\end{lemma}
\begin{proof}
By \reff{nice P} and Lemma~\ref{lemmaA1} in the appendix we have
that
\begin{equation}
\label{multipleM-} %
M_n^-:=\prod_{d=1}^{n-1} \Phi_d(\P{})=\prod_{\substack{d|\R k
\mbox{\scriptsize{ for}}
\\\mbox{\scriptsize{some} } 1\le k\le n-1}}\Phi_d(p)
\end{equation}
is a common multiple of all $\P{j}-1$, $j=1,\ldots,n-1$. Next, for
each $1\le l\le \R-1$, $(l,\R )=1$ we define $1\le b_l \le \R -1$
by $b_l\equiv \r /l \mbox{ mod } \R$. Notice that if
$d\in\mathbb{N}$ satisfies $dl=\R k+\r $ for some
$k\in\mathbb{Z}$, then $d\equiv b_l \mbox{ mod } \R$ and $(l,\R
)=1$ since $(\r ,\R )=1$. Hence
\begin{equation}
\label{Mn} %
\mathcal{M}_n:= \prod_{\substack{d|\R k+\r \mbox{\scriptsize{
for}} \\\mbox{\scriptsize{some} } 0\le k\le
n-1}}\Phi_d(p)=\prod_{\substack{l=1\\(l,\R )=1}}^{\R }\prod_{j=0}^
{\left\lfloor{\frac{n-1}{l}-\frac{lb_l-\r }{l\R
}}\right\rfloor}\Phi_{j\R +b_l}(p)
\end{equation}
is a common multiple of all $\p{}\P{k}-1$, $k=0,\ldots,n-1$. So,
$M_n^-\, \mathcal{M}_n$ is a multiple of $e_n^-$. However, there
are some factors of the form $\Phi_{j\R +b_l}(p)$ appearing in
both $\mathcal{M}_n$ and $M_n^-$. Looking at \reff{multipleM-},
since $(\R ,b_l)=1$ this means that $j\R +b_l$ should be a divisor
of a natural number $k\leq n-1$. So, if $n$ is large enough the
factor $\Phi_{j\R +b_l}(p)$ of $\mathcal{M}_n$ is also present in
$M_n^-$ for $j$ from 0 up to $\lfloor\frac{n-1}{\R }\rfloor -1$
$(\leq\left\lfloor{\frac{n-1}{l}-\frac{lb_l-\r }{l\R
}}\right\rfloor)$, meaning that they have the common factor
\begin{equation}
\label{common} %
C_n^-:=\prod_{\substack{l=1\\(l,\R )=1}}^{\R }
\prod_{j=0}^{\lfloor{\frac{n-1}{\R }}\rfloor -1}\Phi_{j\R
+b_l}(p).
\end{equation}

We proved that $M_n^-\, \mathcal{M}_n \, /\, C_n^-$ is a multiple
of $e_n^-$. Now we look at its asymptotic behaviour. Applying
\reff{Azonderb} on \reff{multipleM-} we easily establish
\begin{equation}
\label{AsM-} %
\log_{\P{}} \left[\lim_{n\to \infty} \left(M_n^-\right)^{1/n^2}
\right] = \frac{3}{\pi^2}.
\end{equation}
Next, recall that $(\R{},b_l)=1$. So, by \reff{Ametb} we also get
\begin{align}
\label{AsM} %
& \log_{\P{}} \left[\lim_{n\to \infty}
\left(\mathcal{M}_n\right)^{1/n^2}
\right] = \frac{3}{\pi^2}\prod_{\substack{\varpi|\R\\
\varpi\: {\rm
prime}}}\frac{\varpi^2}{\varpi^2-1}\sum_{\substack{l=1\\ (l,\R
)=1}}^{\R }\frac{1}{l^2}<\frac{1}{2},\\[1ex]
\label{AsC-} %
& \log_{\P{}} \left[\lim_{n\to \infty} \left(C_n^-\right)^{1/n^2}
\right] = \frac{3}{\pi^2}\frac{\phi(\R)}{\R^2}\prod_{\substack{\varpi|\R\\
\varpi\: {\rm prime}}}\frac{\varpi^2}{\varpi^2-1}=\frac{3}{\pi^2}\frac{1}{\R}\prod_{\substack{\varpi|\R\\
\varpi\: {\rm prime}}}\frac{\varpi}{\varpi+1},
\end{align}
where the last equality follows from the well-known fact
\begin{equation}
\label{propPhi}
\frac{\phi(m)}{m}=\prod_{\substack{\varpi|m\\\varpi\: {\rm
prime}}}\frac{\varpi-1}{\varpi}.
\end{equation}
Combining \reff{AsM-}, \reff{AsM} and \reff{AsC-} we then finally
obtain \reff{AsLCM-}.
\end{proof}
\begin{remark}
The common multiple $\mathcal{M}_n$ of all $\p{}\P{k}-1$,
$k=0,\ldots,n-1$ and its asymptotic behaviour were discussed
already in \cite[Lemma~2]{Matala}.
\end{remark}
\begin{lemma}
\label{lemmaLCM+} %
Let $\r ,\R \in\mathbb{N}$ and $(\r ,\R )=1$. Then
\begin{equation}
\label{AsLCM+} %
\lim_{n\to \infty}\left[{\rm lcm}\left\{\left. \P{j}+ 1, \,
\p{}\P{k}-1 \, \right| \, 1\le j\le n-1, \, 0\le k \le n-1
\right\}\right]^{1/n^2}\le \P{\theta^+ (\R )},
\end{equation}
where
\[
\theta^+(\R )=\gamma^+(\R ) - \frac{3}{\pi^2}\prod_{\substack{\varpi|\R\\
\varpi\: {\rm
prime}}}\frac{\varpi^2}{\varpi^2-1}\sum_{\substack{l=1\\ (l,\R
)=1}}^{\R }\frac{1}{l^2}
\]
with $0<\theta^+(\R ) < \frac{4}{\pi^2}+\frac{1}{2}$ and
$\gamma^+(\R )$ defined as in \reff{gamma+}.
\end{lemma}
\begin{proof}
Recall from the proof of Lemma~\ref{lemmaLCM-} that
$\mathcal{M}_n$, defined as in \reff{Mn}, is a common multiple of
all $\p{}\P{k}-1$, $k=0,\ldots,n-1$. By the property
$x^j+1=(x^{2j}-1)/(x^j-1)$, \reff{nice P} and Lemma~\ref{lemmaA2}
in the appendix we also have that
\begin{equation}
\label{multipleM+} %
M_n^+:=\prod_{d=1}^{n-1} \Phi_{2d}(\P{})=\prod_{\substack{d|2\R
k,\, d \nmid \, \R k \mbox{\scriptsize{ for}}
\\\mbox{\scriptsize{some} } 1\le k\le n-1}}\Phi_d(p)
\end{equation}
is a common multiple of all $\P{j}+1$, $j=1,\ldots,n-1$. So,
$M_n^+\, \mathcal{M}_n$ is a multiple of $e_n^+$.

If $\R{}$ is even, then the index $j\R{}+b_l$ is odd since
$(\R{},b_l)=1$. So, in this case there are no factors of the form
$\Phi_{j\R +b_l}(p)$ appearing in both $\mathcal{M}_n$ and
$M_n^+$. Now suppose that $\R{}$ is odd. Then
\[
M_n^+=\prod_{\substack{d|\R k \mbox{\scriptsize{ for}}
\\\mbox{\scriptsize{some} } 1\le k\le n-1}}\Phi_{2d}(p).
\]
This implies that for a factor $\Phi_{j\R +b_l}(p)$ of
$\mathcal{M}_n$ also appearing in $M_n^+$, we should have $j\equiv
b_l \mod 2$ and $\frac{jr_2+b_l}{2}$ should be a divisor of a
natural number $k\leq n-1$. So, in the case that $\R{}$ is odd,
$\mathcal{M}_n$ and $M_n^+$ have the common factor
\begin{equation} %
\label{Cn+} C_n^+:=\prod_{\substack{l=1\\(l,\R )=1}}^{\R }
\prod_{\substack{j=0\\j\equiv b_l \mbox{\,{\scriptsize mod}\,}
2}}^{\min\left(\lfloor \frac{2(n-1)}{\R }\rfloor -1
,\left\lfloor{\frac{n-1}{l}-\frac{lb_l-\r }{l\R }}\right\rfloor
\right)}\Phi_{j\R +b_l}(p).
\end{equation}
and $M_n^+\, \mathcal{M}_n \, /\, C_n^+$ is a multiple of $e_n^+$.

Now we are interested in the asymptotic behaviour of $M_n^+$ and
$C_n^+$. Applying \reff{Azonderb} on \reff{multipleM+} we easily
obtain
\begin{equation}
\label{AsM+} %
\log_{\P{}} \left[\lim_{n\to \infty} \left(M_n^+\right)^{1/n^2}
\right] = \frac{4}{\pi^2}<\frac{1}{2}.
\end{equation}
Next, we suppose $\R{}$ is odd and look at \reff{Cn+}. Note that
for $n$ large enough we have $\lfloor \frac{2(n-1)}{\R }\rfloor -1
> \left\lfloor{\frac{n-1}{l}-\frac{lb_l-\r }{l\R
}}\right\rfloor$ if and only if $l> \frac{\R}{2}$. Moreover, if
$b_l$ is even, then $j=2i$ is even and we write $j\R +b_l=2(i\R
+b_l/2)$. On the other hand, if $b_l$ is odd then $j=2i+1$ is odd
and we write $j\R +b_l=2(i\R +(b_l+\R{})/2)$. By \reff{A2metb} and
\reff{propPhi} we then get
\begin{align}
\nonumber \log_{\P{}} \left[\lim_{n\to \infty}
\left(C_n^+\right)^{1/n^2}
\right] & = \frac{\phi(\R)}{2\R^2}\frac{4}{\pi^2}\prod_{\substack{\varpi|\R\\
\varpi\: {\rm prime}}}\frac{\varpi^2}{\varpi^2-1}
+\frac{1}{\pi^2}\prod_{\substack{\varpi|\R\\
\varpi\: {\rm prime}}}\frac{\varpi^2}{\varpi^2-1}
\sum_{\substack{l=\lceil \frac{\R{}}{2} \rceil \\(l,\R )=1}}^{\R }
\frac{1}{l^2}\\
\label{AsC+} %
& =\frac{1}{\R}\frac{2}{\pi^2}\prod_{\substack{\varpi|\R\\
\varpi\: {\rm prime}}}\frac{\varpi}{\varpi+1}
+\frac{1}{\pi^2}\prod_{\substack{\varpi|\R\\
\varpi\: {\rm prime}}}\frac{\varpi^2}{\varpi^2-1}
\sum_{\substack{l=\lceil \frac{\R{}}{2} \rceil \\(l,\R )=1}}^{\R }
\frac{1}{l^2}.
\end{align}
Combining \reff{AsM+}, \reff{AsM} and \reff{AsC+} we then finally
obtain \reff{AsLCM+}.
\end{proof}
By \reff{nice P} we obtain
\[
\frac{(\p{}\P{n};\P{})_k}{(\p{};\P{})_{k}} =
\frac{\prod_{i=n}^{n+k-1} \prod_{d|\R i+\r} \Phi_d
(p)}{\prod_{i=0}^{k-1} \prod_{d|\R i+\r} \Phi_d (p)}.
\]
Having a closer look at this expression, it is clear that its
denominator is a divisor of $\mathcal{M}_n$, defined as in
\reff{Mn}, for each $0\le k\le n-1$. As a corollary of
Lemma~\ref{lemmaLCM-}, Lemma~\ref{lemmaLCM+} and \reff{AsM} we
then get the following asymptotic behaviour for $e_n^\pm$.
\begin{corollary}
\label{en asymp} %
Let $\r ,\R \in\mathbb{N}$ and $(\r ,\R )=1$. For $e_n^\pm$
defined as in \reff{en} we have
\begin{equation}
\label{en asymp formule} %
\lim_{n\rightarrow\infty}\left|e_n^\pm\right|^{1/n^2} \le
\P{\gamma^\pm(\R )},
\end{equation}
where $\gamma^-(\R )$ and $\gamma^+(\R )$ are defined as in
\reff{gamma-} and \reff{gamma+}.  As a result we now have
\begin{equation} %
\label{errorasympspeciaal} \lim_{n\rightarrow\infty}\left|b_n^\pm
h^{\pm}-a_n^\pm \right|^{1/n^2}\le \P{\gamma^\pm(\R
)-\frac{3}{2}}.
\end{equation}
\end{corollary}
\subsection{Proof of Theorem~\ref{irrationality}}

In the previous sections we defined integer sequences $a_n^\pm$
and $b_n^\pm$ and managed to find the asymptotic behaviour of
$b_n^\pm h^{\pm}-a_n^\pm$.  Putting these results together, we can
now prove Theorem~\ref{irrationality}.
\begin{varproof}{\bf of Theorem~\ref{irrationality}.}
In Lemma~\ref{lemma en} we made sure that $a_n^\pm$ and $b_n^\pm$,
defined as in \reff{an} and \reff{bn}, are integer sequences.
Note that by \reff{kappa}, \reff{hat Q} and \reff{en asymp
formule} we then get
\begin{equation}
\label{bn asymp lemma nieuw speciaal} \lim_{n\rightarrow\infty}
|b_n|^{1/n^2}\leq \P{\gamma^\pm(\R )+\frac{3}{2}}.
\end{equation}
Lemma~\ref{restterm asymp lemma} assures us that $b_n^\pm
h^{\pm}-a_n^\pm\neq 0$ for all $n\in \mathbb{N}$ and since
$\gamma^\pm<\frac32$, \reff{errorasympspeciaal} guarantees that
$\lim_{n\rightarrow\infty}\left|b_n^\pm h^{\pm}-a_n^\pm\right|=0$.
So, all the conditions of Lemma~\ref{irrlemma} are fulfilled and
$h^{\pm}$ is irrational. Obviously, the irrationality also follows
from the result in the general case, as given in
Theorem~\ref{irrationality alg}. The sharper asymptotics for this
case were only necessary to find the upper bounds for the
irrationality measure, as proposed in Corollary~\ref{irrationality
measure}.
\end{varproof}

\section{Final remark on the special case \reff{specgeval}}
\label{finalremark}

In Lemma~\ref{lemma en} we proposed a possible factor $\e_n^\pm$
such that $a_{n}^\pm$ and $b_{n}^\pm$, defined as in \reff{an} and
\reff{bn}, are integers. It seems that this choice is
(asymptotically) not the optimal one. Empirically (using Maple) we
observed  that $e_{n}^\pm$ can be replaced by
\begin{align}
 \label{xi-}\xi_{n}^- & = \frac{(\p{};\P{})_n}{(\P{};\P{})_{n-1}}\,{\rm lcm}\left\{\left.
\P{j}- 1, \, \p{}\P{k}-1 \,
\right| \, 1\le j\le n-1, \, 0\le k \le n-1 \right\},\\
\nonumber  \xi_{n}^+ & =
\frac{(\p{};\P{})_n}{(\P{};\P{})_{n}}\,\prod_{i=1}^{n}\frac{\P{m_i}-1}{p^{m_i}-1}
\ \prod_{\substack{d|\R,\, d\not=2\\ d \mbox{
\scriptsize{prime}}}} \, \prod_{k=1}^\infty
\Bigl(\phi_{d^k}(p)\Bigr)^{\left\lfloor
\frac{n}{d^k} \right\rfloor} \\ %
\label{xi+}& \qquad \qquad \qquad \qquad \times{\rm
lcm}\left\{\left. \P{j}+ 1, \, \p{}\P{k}-1 \, \right| \, 1\le j\le
n-1, \, 0\le k \le n-1 \right\},
\end{align}
where $m_i$ is the highest odd factor of $i$ and $\lfloor
x\rfloor$ is the integer part of $x$. Then the $a_{n}^\pm$ and
$b_{n}^\pm$ are still integers and ${\rm gcd}
(a_{n}^\pm,b_{n}^\pm)$ turns out to be very small and
asymptotically irrelevant. Up to now we do not have an exact proof
for this, except for the case $\R =1$.  Having in mind the proof
of Lemma~\ref{lemma en}, in this case the denominator of
\[\frac{(\p{}p^n;p)_k}{(\p{};p)_k}=\left[{n+k+r_1-1\atop n}\right]_{p}\left[{n+r_1-1\atop n}\right]_{p}^{-1}=\left[{n+k+r_1-1\atop n}\right]_{p}\frac{(p;p)_n}{(\p{};p)_{n}}\]
is clearly cancelled out by the first factor of \reff{xi-} and
\reff{xi+}.

Now define $\delta^-(\R ):=\theta^-(\R )$ and $\delta^+(\R
):=\theta^+(\R )+\frac{1}{3}-\frac{1}{3\R}$, where $\theta^\pm(\R
)$ is as in Lemma~\ref{lemmaLCM-} and Lemma~\ref{lemmaLCM+}.
Following the arguments of this paper, with the adjustment
mentioned above we could prove the existence of integer sequences
$\alpha_n^\pm$, $\beta_n^\pm$ such that
\begin{align*}
& \lim_{n\rightarrow\infty}\left|\xi_n^\pm\right|^{1/n^2} \le
\P{\delta^\pm(\R )}, \\
& \left|\beta_{n}^\pm h^\pm -\alpha_{n}^\pm\right|
=\O\Bigl((\beta_n^\pm)^{-\frac{3-2\,\delta^\pm(\R
)}{3+2\,\delta^\pm(\R )}+\varepsilon}\Bigr),\qquad \mbox{for all }
\varepsilon
>0,\qquad n\to \infty,
\end{align*}
implying
\begin{equation}
\label{upbfinal} \mu\left(h^\pm\right)\le
\chi^\pm(\R):=\left(\frac{3-2\,\delta^\pm(\R )}{6}\right)^{-1}.
\end{equation}
Comparing Table~\ref{tabel1} with Table~\ref{tabel2} we see that
this would considerably improve the upper bound for the
irrationality measure.

\begin{table}[t]
\begin{center}
\begin{tabular}{|c|rcl|rcl|}
\hline %
$\R$ &  \multicolumn{3}{|c|}{$\chi^-(\R)$}&\multicolumn{3}{|c|}{$\chi^+(\R)$}\\
\hline & & & & & & \\[-2ex]
1 & $\frac{2\pi^2}{\pi^2-2}$&$\approx$&$ 2.508284762
$&$\frac{6\pi^2}{3\pi^2-8}$&$\approx$&$ 2.740438628 $ \\[0.5ex]
2 & $\frac{2\pi^2}{\pi^2-4}$&$\approx$&$ 3.362953864
$&$\frac{9\pi^2}{4\pi^2-24}$&$\approx$&$ 5.738728718 $ \\[0.5ex]
3 & $\frac{32\pi^2}{16\pi^2-69}$&$\approx$&$ 3.552067296
$&$\frac{432\pi^2}{184\pi^2-1071}$&$\approx$&$ 5.722990389 $ \\[0.5ex]
4 & $\frac{54\pi^2}{27\pi^2-125}$&$\approx$& $3.767042717
$&$\frac{108\pi^2}{45\pi^2-304}$&$\approx$&$ 7.606512209 $ \\[0.5ex]
5 & $\frac{3456\pi^2}{1728\pi^2-8005}$&$\approx$&$ 3.769124031
$&$\frac{25920\pi^2}{10656\pi^2-68555}$&$\approx$&$ 6.986661787 $ \\[0.5ex]
6 & $\frac{300\pi^2}{150\pi^2-743}$&$\approx$& $4.015077389
$&$\frac{675\pi^2}{275\pi^2-1953}$&$\approx$&$ 8.752624192 $ \\[0.5ex]
7 & $\frac{172800\pi^2}{86400\pi^2-414281}$&$\approx$&$
3.889740382
$&$\frac{1814400\pi^2}{734400\pi^2-4949917}$&$\approx$&$
7.791520131 $
\\[0.5ex]
8 & $\frac{132300\pi^2}{66150\pi^2-327931}$&$\approx$&
$4.018388861
$&$\frac{264600\pi^2}{106575\pi^2-766112}$&$\approx$&$ 9.139383255 $ \\[0.5ex]
9 & $\frac{1881600\pi^2}{940800\pi^2-4664047}$&$\approx$&$
4.018510114
$&$\frac{25401600\pi^2}{10192000\pi^2-71413173}$&$\approx$&$
8.592266790 $
\\[0.5ex]
10 & $\frac{71442\pi^2}{35721\pi^2-180973}$&$\approx$&
$4.109498873
$&$\frac{178605\pi^2}{71442\pi^2-521890}$&$\approx$&$ 9.621306357 $ \\[0.5ex]
\hline%
\end{tabular}
\caption{\label{tabel2} Some values of $\chi^\pm(\R)$, see
\reff{upbfinal}.}
\end{center}
\end{table}

\section*{Acknowledgements}

We like to thank W. Van Assche for careful reading and useful
discussions.

This work was supported by INTAS project 03-51-6637, by FWO
projects G.0184.02 and G.0455.05 and by OT/04/21 of K.U.Leuven.
The first author is a postdoctoral researcher at the K.U.Leuven
(Belgium).

\appendix

\section*{Appendix}

\setcounter{section}{1}
\setcounter{equation}{0}

In this appendix we prove some properties dealing with the
cyclotomic polynomials \reff{cyclotomic def}. In the proofs we use
the notation $a|_c^{} b$ which means that $a|b$ and that $a$
contains any prime factor present in both $b$ and $c$, up to the
highest possible power.
\begin{lemma}
\label{lemmaA1}
The cyclotomic polynomials satisfy
\begin{equation}
\label{app gelijk1}
\prod_{d=1}^{n-1}\Phi_d(x^r)=\prod_{\substack{d|rk \:{\rm
for}\\{\rm some}\: 1\leq k\leq n-1}}\Phi_d(x), \qquad
n\in\mathbb{N}.
\end{equation}
\end{lemma}
\begin{proof}
Both expressions in \reff{app gelijk1} are a common multiple of
the polynomials $\{x^{rj}-1,j=1,\ldots,n-1\}$, the latter one
being the least common multiple.  Since they are both monic
polynomials, it is sufficient to prove that their respective
degrees are equal.  Recall that Euler's totient function $\phi$
represents the degree of the cyclotomic polynomials. So, we have
to prove the equality
\[r\sum_{d=1}^{n-1}\phi(d)=\sum_{\substack{d|rk \:\:{\rm for}\\{\rm some}\: 1\leq k\leq
n-1}}\phi(d).\]
It is easily seen that this holds for every $n\in\mathbb{N}$ if
and only if
\begin{equation}
\label{totientTB1} %
r\phi(n)=\sum_{d|_n^{} r}\phi(nd),\qquad n\in\mathbb{N}.
\end{equation}
Denote by $g$ the smallest divisor of $r$ such that $r/g$ contains
no prime factors present in $n$. Next, recall some well-known
properties of the totient function, such as
\begin{align}
\label{eig1}
   & \phi(ab)=\phi(a)\phi(b) \qquad \hspace{.485cm}\mbox{ if } (a,b)=1;\\
\label{eig2}
   & \phi(ab)=a\phi(b) \qquad \hspace{1cm}\mbox{ if $a$ prime and $a|b$}.
\end{align}
Using these and \reff{nice P} we obtain
\[\sum_{d|_n^{} r}\phi(nd)=\sum_{d|\frac{r}{g}}\phi(ndg)
=\phi(ng)\sum_{d|\frac{r}{g}}\phi(d)
=\phi(ng)\frac{r}{g}=r\phi(n),\]
which is \reff{totientTB1}. This then proves the statement of the
lemma.
\end{proof}

\begin{lemma}
\label{lemmaA2}
For any prime number $\sigma$ we have
\begin{equation}
\label{app gelijk2} %
\prod_{d=1}^{n-1}\Phi_{\sigma d}(x^r)=\prod_{\substack{d|\sigma
rk,\, d\nmid rk \:{\rm for}\\{\rm some}\: 1\leq k \leq
n-1}}\Phi_d(x),\qquad n\in\mathbb{N}.
\end{equation}
\end{lemma}
\begin{proof}
First of all note that both expressions are a common multiple of
the set of polynomials
\[\left\{\frac{x^{\sigma rj}-1}{x^{rj}-1},j=1,\ldots,n-1\right\},\]
the latter one being the least common multiple.  As in the proof
of the previous lemma, since they are both monic polynomials, it
is sufficient to prove that they have the same degree. Hence we
have to prove the equality
\[r\sum_{d=1}^{n-1}\phi(\sigma
d)=\sum_{\substack{d|\sigma rk,\,d\nmid rk \:{\rm for}\\{\rm
some}\: 1\leq k \leq n-1}}\phi(d).\]
Introduce the notation $r=\sigma^\tau r'$ where $r'$ does not
contain the prime factor $\sigma$. This then holds for every
$n\in\mathbb{N}$ if and only if
\begin{equation}
\label{totientTB2} %
r\phi(\sigma n)=\sum_{d|_{n}^{} r'}\phi(\sigma^{\tau+1} nd),\qquad
n\in\mathbb{N}.
\end{equation}
From this point, we can proceed in an analogous way as in the
previous lemma. Denote by $g$ the smallest divisor of $r'$ such
that $r'/g$ contains no prime factors present in $n$. By
\reff{eig1} we then get
\[
\sum_{d|_{n}^{} r'}\phi(\sigma^{\tau+1}
nd)=\sum_{d|\frac{r'}{g}}\phi(\sigma^{\tau+1}
ngd)=\phi(\sigma^{\tau+1} ng)\sum_{d|\frac{r'}{g}}\phi(d).
\]
Next, applying \reff{nice P} and \reff{eig2} we finally obtain
\[
\sum_{d|_{n}^{} r'}\phi(\sigma^{\tau+1} nd) =\phi(\sigma^{\tau+1}
ng)\frac{r'}{g}=\phi(\sigma^{\tau+1} n)r'=\phi(\sigma n)r.
\]
This ends the proof of \reff{totientTB2} and hence of this lemma.
\end{proof}

\end{document}